\documentclass[12pt]{article} 
\usepackage{amsmath, amscd, amssymb,latexsym, epsfig, color, amsthm, enumerate}

\newtheorem{theorem}{Theorem}[section]

\newtheorem{corollary}[theorem]{Corollary}

\newtheorem{lemma}[theorem]{Lemma}

\newtheorem{example}[theorem]{Example}
\newtheorem{question}[theorem]{Question}

\theoremstyle{definition}

\theoremstyle{remark}

\def\proof{\smallskip\noindent {\it Proof: \ }} 
\def\endproof{\hfill$\square$\medskip}
 \def\Z{\mathbb{Z}} 
 
\def\Q{\mathbb{Q}}

\def\a{{\bf{a}}}

\def\n{{\bf{n}}}
 
\def\X{\mathbf{X}}
\def\Y{\mathbf{Y}}
\def\V{\mathbf{V}}
\def\Res{\mathcal{R}}
\def\I{\mathcal{I}}

\DeclareMathOperator{\Lex}{Lex}
\DeclareMathOperator{\Init}{In}
\DeclareMathOperator{\Skel}{Skel}

\newcommand{\field}{{\bf k}} 
\newcommand{\Field}{{\bf K}}

\newcommand{\Char}{\mbox{\upshape char}\,}

\newcommand{\mult}[3]{\operatorname{S}_{#3}(#1, (\infty, #2))}
\newcommand{\minv}{\mathfrak{m}}
\newcommand{\comp}[3]{\Lambda_{#3}(#1, #2)}

\DeclareMathOperator{\fac}{fac}
\DeclareMathOperator{\FL}{full}
\DeclareMathOperator{\fgap}{fgap}
\DeclareMathOperator{\up}{low}
\DeclareMathOperator{\Gap}{Gap}
\DeclareMathOperator{\tail}{tail}
\DeclareMathOperator{\rlex}{rl}
\newcommand{\fll}[3]{{#1}[#2; #3]}
\newcommand{\bk}[2]{#1\langle #2  \rangle}

\title{Face numbers of generalized balanced Cohen--Macaulay complexes}

\author{Jonathan Browder and Isabella Novik
\thanks{Novik's research is partially supported by 
Alfred P.~Sloan Research Fellowship and NSF grant DMS-0801152}\\ 
\small Department of Mathematics, Box 354350\\[-0.8ex] 
\small University of Washington, Seattle, WA 98195-4350, USA,\\[-0.8ex] 
\small \texttt{[browder, novik]@math.washington.edu} }

\begin{document}

\maketitle 

\begin{abstract} 
A common generalization of two theorems on 
the face numbers of Cohen--Macaulay (CM, for short) 
simplicial complexes is established: 
the first is the theorem of Stanley (necessity)
and Bj\"orner-Frankl-Stanley (sufficiency) that characterizes 
all possible face numbers of $\a$-balanced CM
complexes, while the second is the theorem of Novik (necessity)
and Browder (sufficiency) that characterizes the face numbers of 
CM subcomplexes of the join of the boundaries of simplices.
\end{abstract}

\section{Introduction}
A basic invariant of a simplicial complex $\Delta$
is its $f$-vector, 
$f(\Delta)=(f_{-1}, f_0, \ldots, f_{\dim \Delta})$,
where $f_i$ denotes the number of $i$-dimensional faces of $\Delta$.
Can one characterize the set of all $f$-vectors of various
interesting families of simplicial complexes? 

In the mid-sixties,
Kruskal \cite{Krus} and Katona \cite{Kat} (independently)
provided an answer for the family of all simplicial complexes.
Their result started the still continuing 
quest for finding characterizations of
the $f$-vectors of other important subfamilies of complexes. 
In particular, Stanley \cite{St77} characterized the 
$f$-vectors of all Cohen--Macaulay (CM, for short) simplicial 
complexes. His theorem was then refined 
in  \cite{St79} (necessity) and \cite{BjFrSt}
(sufficiency) to a characterization of all possible $f$-vectors
of $\a$-balanced CM complexes. 
More recently, in connection to complexes
endowed with a proper group action, 
the class of CM subcomplexes of the join of the
boundaries of simplices was considered,
and its collection of $f$-vectors characterized --- 
see \cite{Nov05} (necessity) and \cite{Br} (sufficiency).
Our goal here is to provide a simultaneous and natural
generalization of both of these results. 
We thank Anders Bj\"orner for suggesting 
that such a generalization might exist.
 
To state our main result we need a bit of notation. 
The complexes we consider are full-dimensional 
subcomplexes of
\begin{equation} \label{Lambda}
\Lambda_d=\Lambda_d(\n, \a)=\Lambda_d(\V_1, \ldots, \V_m, \a):=
\Skel_{d-1}\left(\Skel_{a_1-1}\overline\V_1
\ast\cdots\ast \Skel_{a_m-1}\overline\V_m\right).
\end{equation}
Here $\n:=(n_1,\ldots, n_m)\in\Z_+^m$ and 
$\a:=(a_1, \ldots, a_m) \in \Z_+^m$
are $m$-dimensional integer vectors,  
$\V_1, \ldots, \V_m$ are pairwise disjoint
sets of sizes $n_1,\ldots, n_m$, resp.,
$\overline\V_i$ denotes the simplex on the vertex set $\V_i$, 
$\ast$ stands for the join of simplicial complexes, and
$\Skel_{a-1}(\Delta)$ is the $(a-1)$-dimensional 
skeleton of $\Delta$ --- the subcomplex of $\Delta$ 
consisting of all faces of $\Delta$ whose dimension 
is strictly smaller than $a$. 

Modulo the results of \cite{BjFrSt, St77, St79}, it is perhaps not
so surprising that
our characterization depends on the notion of a multicomplex.
Specifically, let $\X_0, \X_1, \ldots, \X_m$ be pairwise disjoint
sets of variables,  $\X$  their union, and  $a_0=\infty$. 
For $0\leq i \leq m$, denote by 
$S(\X_i, a_i)$ the poset of all monomials on $\X_i$ (ordered by
divisibility) of degree at most $a_i$. Let $S(\X, (\infty, \a))$ be 
the product of these posets, and  
$S=S_d(\X, (\infty, \a))$ 
its subset of monomials of degree at most $d$:
that is, $S$ consists of all monomials 
$\mu=\mu_{\X_0}\mu_{\X_1}\cdots\mu_{\X_m}$ of degree
no greater than $d$, where for $0\leq i \leq m$, $\mu_{\X_i}$ 
is a monomial on $\X_i$ of degree no greater than $a_i$. For $\mu\in S$ and
$Y\subseteq \X$, we also denote by $\mu_Y$ the part of $\mu$ 
supported in $Y$.

The main result of this paper is the following.
We postpone the discussion of the relevant definitions, 
including those of the $h$- and $F$-vectors, to the next section.
\begin{theorem} \label{main}
Let $d$ be a positive integer and let 
$h=(h_0, \ldots, h_{d})\in\Z^{d+1}$. 
If $\a=(a_1, \ldots, a_m)\in\Z^m$ and 
$\n=(n_1,\ldots, n_m)\in\Z^m$ are such that
$n_i\geq a_i>0$ for $1\leq i \leq m$ and 
$\sum_{i=1}^m a_i \geq d$, then the following are equivalent:

\begin{enumerate}
\item $h$ is the $h$-vector of a $(d-1)$-dimensional 
$\Q$-CM subcomplex of $\Lambda_d=\Lambda_d(\n, \a)$,

\item $h$ is the $h$-vector of a $(d-1)$-dimensional 
shellable subcomplex of $\Lambda_d=\Lambda_d(\n, \a)$,

\item $h$ is the $F$-vector of a multicomplex 
$M\subseteq S_d(\X, (\infty, \a))$, where 
$|\X_0|=\left(\sum_{i=1}^m a_i\right) -d$ and 
$|\X_i|=n_i-a_i$ for all $1\leq i \leq m$,

\item $h$ is the $F$-vector of a $(0)$-compressed multicomplex 
$M\subseteq S_d(\X, (\infty, \a))$, where $|\X_0|=\left(\sum_{i=1}^m a_i\right) -d$
and $|\X_i|=n_i-a_i$ for all $1\leq i \leq m$.
\end{enumerate}
\end{theorem}

\noindent
In the case of $\sum_{i=1}^m a_i = d$, the full-dimensional
subcomplexes of $\Lambda_d$ 
are precisely the $\a$-balanced complexes of \cite{BjFrSt, St79} 
and Theorem \ref{main} reduces to (the $f$-vector rather 
than flag $f$-vector version of) \cite[Theorem 1]{BjFrSt},
while in the case of $n_i-a_i\in\{0,1\}$ for all $1\leq i \leq m$,
one obtains complexes considered in \cite{Br, Nov05} and
recovers \cite[Cor.~1]{Br}. 

A pure simplicial complex is completely balanced if it is
$\a$-balanced with $\a=(1,1,\ldots, 1)$.
It is worth remarking that for completely balanced CM
complexes, the main result of \cite{FrFuKa} turns 
Theorem \ref{main} into a numerical characterization 
of the $h$-numbers of such complexes. Similarly, in the case of
$n_i-a_i\in\{0,1\}$ for all $1\leq i \leq m$, the Clements-Lindstr\"om
theorem \cite{ClLi} combined with Theorem \ref{main} provides a
numerical characterization of the $h$-numbers of such CM complexes.
However, for a general $\a$,
no numerical characterization of the $F$-vectors of sub-multicomplexes
of $S_d(\X, (\infty,\a))$ is known at present. 

\begin{question} Can one use the combinatorial characterization of the
h-numbers given by Theorem \ref{main} to arrive at a numerical 
characterization of the $h$-numbers?
\end{question}

The rest of the paper is structured as follows. Section 2 contains
basics on simplicial complexes and Stanley-Reisner rings. 
The implication $1.\rightarrow 3.$ is proved in Section 3. 
The main technique employed in the proof of this part is the study of
combinatorics of the (non-generic) initial ideal of the Stanley-Reisner 
ideal of the complex in question.
The notion of $(0)$-compressed multicomplexes is introduced in Section 4.
The same section contains the proof of $3.\rightarrow 4.$ part.
The implication $4.\rightarrow 2.$ is the most technical one and
is established in Sections 5 and 6. Finally, the implication
$2.\rightarrow 1.$ follows from the well-known fact 
\cite[Thm.~III.2.5]{St96} that every shellable complex is CM. The flavor
of the proofs of $3.\rightarrow 4.$ and $4.\rightarrow 2.$ is 
motivated by and somewhat similar to the proofs in \cite{BjFrSt}.

\section{Preliminaries}
Here we briefly review several notions and results
related to simplicial complexes and Stanley-Reisner rings.
An excellent reference to this material is 
Stanley's book \cite{St96}.

\smallskip\noindent{\bf Complexes and multicomplexes.}
A {\em simplicial complex} $\Delta$ on the vertex set $\V$ is 
a collection of subsets of $\V$ that is closed under inclusion.
(We do not require that $\Delta$ contains all singletons 
$\{v\}$ for $v\in \V$.) 
The elements of $\Delta$ are called its {\em faces}.
For $\tau\in\Delta$, set $\dim \tau:=|\tau|-1$ and define the 
{\em dimension} of $\Delta$, $\dim \Delta$, as the maximal 
dimension of its faces. The {\em facets} of $\Delta$
are maximal (under inclusion) faces of $\Delta$. 
We say that $\Delta$ is {\em pure} if all of its facets have
the same dimension. The $f$-vector of $\Delta$ is  
$f(\Delta)=(f_{-1}, f_0, \ldots, f_{d-1})$, 
where $d-1=\dim\Delta$ and $f_j$ is the number of $j$-dimensional 
faces of $\Delta$. It is sometimes more convenient to work 
with the $h$-vector, $h(\Delta)=(h_0,h_1,\ldots, h_d)$, that
carries the same information as $f(\Delta)$ and is defined by
$\sum_{i=0}^d h_i x^{d-i}= \sum_{i=0}^d f_{i-1} (x-1)^{d-i}$.

Similarly, a {\em multicomplex} $M$ on the set of variables $\X$ is
a collection of monomials supported in $\X$ that is closed under
divisibility. 
The $F$-vector of
a multicomplex $M$ is the vector $F(M)=(F_0, F_1, \ldots)$,
where $F_j=F_j(M):=|\{\mu\in M \ : \ \deg\mu=j\}|$.
Thus a simplicial complex can be naturally
identified with a multicomplex all of whose elements 
are squarefree monomials. Under this identification, 
the $f$-vector of a simplicial complex differs from its $F$-vector
only by a shift in the indexing.

If $\Delta_1$ and $\Delta_2$ are simplicial complexes on disjoint
vertex sets $\V_1$ and $\V_2$, then their {\em join} is the following
simplicial complex on $\V_1\cup \V_2$: 
$\Delta_1\ast\Delta_2:= 
\{\tau_1\cup \tau_2 \ : \ \tau_1\in\Delta_1, \ \tau_2\in \Delta_2\}.$

\medskip\noindent{\bf Shellability.}
Let $\Delta$ be a pure $(d-1)$-dimensional
simplicial complex. For $\tau\in\Delta$, denote
by $\overline{\tau}$ the simplex $\tau$ together with all its faces.
A {\em shelling} of $\Delta$ is an ordering 
$(\tau_1, \tau_2, \ldots, \tau_s)$ 
of its facets such that for all $1<i\leq s$, the complex
$\overline{\tau}_i\cap(\cup_{j<i}\overline{\tau}_j)$
is pure of dimension $d-2$. Such an ordering is then called a
{\em shelling}. Equivalently, $L=(\tau_1, \tau_2, \ldots, \tau_s)$ is
a shelling if for every $1\leq i\leq s$, the face 
$\Res_L(\tau_i):=\{v\in\tau_i \ : \ \tau_i-\{v\} 
\subseteq \tau_j \mbox{ for some } j<i \}$ is the unique minimal
face of $\overline{\tau_i}-(\cup _{j<i}\overline{\tau}_j)$, called 
the {\em restriction} of $\tau_i$.
It was realized by McMullen \cite{Mc} that 
the $h$-vector of a shellable complex $\Delta$ can be easily computed 
from its shelling $L=(\tau_1, \tau_2, \ldots, \tau_s)$: 
\begin{equation}  \label{h-shellable}
h_i = \left|\{\tau_j \ : \ |\Res_L(\tau_j)|=i\}\right|, \quad i=0,1,\ldots, d.
\end{equation}

\medskip\noindent
{\bf Stanley-Reisner rings and Cohen-Macaulay complexes.} 
Let $\Delta$ be a simplicial complex on 
the vertex set $\V$, and let $\tilde{\X}=\{x_v \, : \, v\in\V\}$
be the corresponding set of variables. 
Fix a field $\field$ of characteristic 0 (e.g. $\Q$), and 
consider $\field[\tilde{\X}]$ --- the polynomial ring over $\field$ 
in variables $\tilde{\X}$ with the grading $\deg x=1$ for $x\in\tilde\X$.
The {\em Stanley-Reisner ideal} of $\Delta$, $I_\Delta$, is the ideal
generated by the squarefree monomials corresponding to non-faces:
$$I_\Delta=(x_{v_1}\cdots x_{v_j} \ : \ \{v_1,\ldots, v_j\}\notin\Delta).$$
The {\em Stanley-Reisner ring} (also known as the 
{\em face ring}) of $\Delta$ is $\field[\Delta]:=\field[\tilde\X]/I_\Delta$.
This ring is graded, and we denote by $\field[\Delta]_i$ its $i$th
homogeneous component.

If $\Delta$ is a $(d-1)$-dimensional simplicial complex, 
then by \cite[Theorem II.1.3]{St96}, the {\em Krull
dimension} of $\field[\Delta]$ (i.e., the maximum number of
algebraically independent elements over $\field$ in $\field[\Delta]$)
is $d$. Moreover, since $\field$ is infinite and $\field[\Delta]$ 
is generated (as an algebra) by $\field[\Delta]_1$, it follows
from the Noether Normalization Lemma that there exists a sequence
$\theta_1, \ldots, \theta_d \in \field[\Delta]_1$  of $d$ 
linear forms such that $\field[\Delta]/(\theta_1, \ldots, \theta_d)$ 
is a finite-dimensional $\field$-space. Such a sequence is called a 
{\em linear system of parameters} (l.s.o.p.) for $\field[\Delta]$.

A simplicial $(d-1)$-dimensional complex $\Delta$ is 
{\em $\field$-Cohen-Macaulay} 
(CM, for short) if for every l.s.o.p.~$\theta_1, \ldots, \theta_d$, 
one has
\begin{equation} \label{h-CM}
\dim_\field \left(\field[\Delta]/(\theta_1, \ldots, \theta_d)\right)_i
=h_i(\Delta) \quad \mbox{for all } 0\leq i \leq d.
\end{equation}
We refer our readers to Chapter II of \cite{St96}
for several other equivalent definitions of CM complexes.
One of them is a result of Reisner \cite{Reis} 
that Cohen-Macaulayness of a complex is equivalent to 
vanishing of certain simplicial homologies. Some immediate
corollaries of Reisner's result are
(i) a CM complex is pure, and (ii)
any triangulated sphere is CM. In addition,
one can use Reisner's criterion to show that all shellable complexes
are CM.

\medskip\noindent{\bf Revlex order and initial ideals.} 
In several parts of the proof of 
Theorem \ref{main} we use the notion of the 
(homogeneous) {\em reverse lexicographic order} 
(on sets and on monomials), which we review now. 
Fix a total order $\succ$ on $\V$. For $S,T\subseteq \V$, 
we 
write $S\succ_{\rlex} T$ (or simply $S\succ T$), if $|S|=|T|$ and 
the least element of $(S-T)\cup(T-S)$ w.r.t.~$\succ$ is in $T$. 
For instance, if  $1\succ 2\succ 3\succ 4$, then 
$\{1,2\}\succ\{1,3\}\succ\{2,3\}\succ\{1,4\}\succ\{2,4\}\succ\{3,4\}$.
Similarly, for a total order $\succ$ on $\tilde\X$, and two monomials
$\mu_1$ and $\mu_2$ on $\tilde\X$, $\mu_1$ is revlex larger
than $\mu_2$ if $\deg\mu_1=\deg\mu_2$ and the least variable 
(w.r.t~$\succ$) that appears in $\mu_1/\mu_2$ has a negative exponent.
Thus, if $x_1\succ x_2\succ x_3$, then 
$x_1^2\succ x_1x_2 \succ x_2^2 \succ x_1x_3\succ x_2x_3\succ x_3^2$.

If $\succ$ is a fixed order on $\tilde\X$ and 
$I$ is a homogeneous ideal of $\field[\tilde\X]$, 
then denote by $\Init(I)$ the
{\em reverse lexicographic initial ideal} of $I$ --- 
the monomial ideal generated by the  
leading (w.r.t~revlex order) monomials 
of all homogeneous elements of $I$ (\cite[Section 15.2]{Eis}).

\section{From CM complexes to multicomplexes}
{\bf Proof outline.}
Given a $(d-1)$-dimensional CM complex $\Delta\subseteq\Lambda_d$, 
how does one construct a multicomplex $M_\Delta\subseteq 
S=S_d(\X, (\infty,\a))$ 
whose $F$-vector equals the $h$-vector of $\Delta$? 
The idea (introduced in \cite{Nov05})
is very simple: instead of working with the 
ideal $I_\Delta$, consider a suitably chosen graded automorphism 
$g: \field[\tilde\X] \rightarrow \field[\tilde\X]$ 
and the ideal $gI_\Delta$ --- the image of $I_\Delta$ under $g$.
Fix a total order $\succ$ on $\tilde\X$, 
and let $L\subset \tilde\X$ be
the set of $d$ last (w.r.t.~$\succ$) variables of $\tilde\X$. 
Then $J_\Delta:=gI_\Delta+(x : x\in L)$ is a homogeneous ideal, 
and hence $\Init(J_\Delta)$ ---  the revlex initial ideal of 
$J_\Delta$ --- is well-defined.

Define $M_\Delta$ to be the set of all monomials on $\tilde\X$ 
that do not belong to $\Init(J_\Delta)$. Then $M_\Delta$ is a 
multicomplex on $\X:=\tilde\X-L$.
This follows immediately from the fact that $\Init(J_\Delta)$ is
an ideal and that all $x\in L$ are elements of $J_\Delta$. 
What are the $F$-numbers of $M_\Delta$?  If $g$ is ``generic enough" 
so that $\{g^{-1}(x) : x\in L\}$ 
is an l.s.o.p.~for $\field[\Delta]=\field[\tilde\X]/I_\Delta$, 
then the elements of $L$ form an l.s.o.p.~for 
$\field[\tilde\X]/gI_\Delta$. 
As $\field[\tilde\X]/gI_\Delta$ is an isomorphic image of 
$\field[\Delta]$, eq.~(\ref{h-CM})  yields that 
$\dim_\field (\field[\tilde\X]/J_\Delta)_i=
h_i(\Delta)$ for all  $0\leq i \leq d$. Theorem 15.3 of \cite{Eis}
asserting that $M_\Delta$ is a $\field$-basis 
for $\field[\tilde\X]/J_\Delta$ then implies that
$F_i(M_\Delta)=\dim_\field (\field[\tilde\X]/J_\Delta)_i=h_i(\Delta)$. 
On the other hand, by choosing $g$ to be not completely generic but in a way that
 ``respects the structure" of  $\Lambda_d$, certain monomials
can be forced to be in $\Init(J_\Delta)$, and hence not in $M_\Delta$, 
ensuring that $M_\Delta$ a subset of $S$.  
Specifying $\succ$ and $g$, and 
verifying that for a given CM subcomplex $\Delta\subseteq\Lambda_d$ 
the above procedure  does produce a multicomplex $M_\Delta\subseteq S$
with $F(M_\Delta)=h(\Delta)$ is the goal of this section.

\medskip\noindent{\bf Specifying $\mathbf\succ$.}
Recall the definition of 
$\Lambda_d=\Lambda_d(\n, \a)=\Lambda_d(\V_1, \ldots, \V_m, \a)$, 
see eq.~(\ref{Lambda}). It is a simplicial complex on 
$\V=\cup_{i=1}^m \V_i$, where $|\V_i|=n_i\geq a_i$. Thus for each
$1\leq i \leq m$, we can split $\V_i$ into two disjoint sets: $V_i$ 
of size $n_i-a_i$ and $W_i$ of size $a_i$. Denote the 
set of variables corresponding to vertices in $V_i$ by $\X_i$
and to vertices in $W_i$ by $\Y_i$, 
$1\leq i \leq m$. Let $\X_{m+1}:= \cup_{i=1}^m \Y_i$
 and $\tilde\X:=\cup_{i=1}^{m+1} \X_i$.
Fix any total order $\succ$ on $\tilde\X$ with the property that for 
$x\in\X_i$ and $y\in \X_j$, $x\succ y$ if $1\leq i<j \leq m+1$.
(In other words, the order $\succ$ first lists the elements of $\X_1$, 
then those of $\X_2, \ldots, \X_m, \X_{m+1}$.)

\medskip\noindent{\bf Defining $g$.}
For this part we replace $\field$ by a larger field $\Field$ ---
the field of rational functions over $\field$ in indeterminates
$\cup_{i=1}^{m+1} \{\alpha^i_{xy}  :  x, y\in \X_i\}\cup 
\cup_{i=1}^m \{\beta^i_{xy}  :  x\in\X_i, y\in \Y_i\}$,
and perform all computations inside $\Field[\tilde\X]$
rather than $\field[\tilde\X]$. In particular, we regard $I_\Delta$
and $I_{\Lambda_d}$ as ideals of $\Field[\tilde\X]$.

Let $A_i=(\alpha^i_{xy})_{x, y\in \X_i}$ ($i=1,\ldots, m$) be an 
$(n_i-a_i)\times (n_i-a_i)$-matrix
whose rows and columns are indexed by elements of $\X_i$, 
let $B_i=(\beta^i_{xy})_{x\in \X_i,y\in\Y_i}$ ($i=1,\ldots, m$) be an
$(n_i-a_i)\times a_i$-matrix
whose rows are indexed by elements of $\X_i$ and columns by $\Y_i$,
and let $C=(\alpha^{m+1}_{xy})_{x,y\in\X_{m+1}}$ be a 
$(\sum_{i=1}^m a_i)\times (\sum_{i=1}^m a_i)$-matrix
whose rows and columns are indexed by elements of $\X_{m+1}$.
Define $A$ to be a block-diagonal square 
matrix that has  $A_1, A_2,\ldots, A_m$ on the main
diagonal and zeros everywhere else. Similarly, let $B$ be 
a block-diagonal (but not square) matrix that has blocks
$B_1, \ldots, B_m$ on the main diagonal and zeros everywhere else.
Then $ABC$ is a well-defined matrix whose rows are
indexed by $\cup_{i=1}^m \X_i$ and whose columns are
indexed by $\X_{m+1}$. This leads to the main definition of this 
section --- a square matrix $g$ whose rows and columns are indexed
by $\tilde\X$. We define
$$
g^{-1}:=\left[\begin{array}{cc}
-A & ABC \\
O & C
\end{array}
\right], \quad
\mbox{so that} \quad
g=\left[\begin{array}{cc}
-A^{-1} & B\\
O & C^{-1}
\end{array}
\right],
$$
where $O$ stands for the zero-matrix.
As $g$ is invertible, it defines a graded automorphism of 
$\Field[\tilde\X]$ via $g(x)=\sum_{y\in\tilde\X} g_{xy}y$.

Let $L$ be the set of least (w.r.t.~$\succ$) $d$ variables of 
$\tilde\X$. Since $|\X_{m+1}|=\sum_{i=1}^m a_i \geq d$, 
$L\subseteq \X_{m+1}$. The following two lemmas
are the main steps in the proof of the implication $1.\rightarrow 3.$
of Theorem \ref{main}.
For subsets $G, H$ of $\tilde\X$, we denote
by $g^{-1}_{G,H}$ the submatrix of $g^{-1}$ whose rows are indexed by $G$ 
and columns by $H$.

\begin{lemma}  \label{KiKl}
If $\tau$ is an arbitrary facet of $\Lambda_d$ 
and $H\subset\tilde\X$ is the corresponding set of variables, then
the submatrix $g^{-1}_{H,L}$ is nonsingular.
\end{lemma}

\begin{lemma} \label{gin}
The ideal $\Init(gI_{\Lambda_d})$ contains all monomials of degree 
$a_i+1$ that are supported in $\X_i$ for all $1\leq i \leq m$.
\end{lemma}

\smallskip\noindent{\it Proof of $1.\rightarrow 3.$: \ }
Using the above notation, let $\X_0=\X_{m+1}-L$, $\X=\cup_{i=0}^m\X_i$,
and define $M_\Delta$ as in the proof outline. Since $\Delta$ is a 
subcomplex of $\Lambda_d$, $I_\Delta\supseteq I_{\Lambda_d}$. Thus,
$J_\Delta \supseteq J_{\Lambda_d} \supseteq gI_{\Lambda_d}$, and hence 
$\Init(J_\Delta) \supseteq \Init(gI_{\Lambda_d})$, which together with
Lemma \ref{gin} implies that for $i=1, \ldots, m$, no monomial of degree 
$a_i+1$ that is supported in $\X_i$ belongs to $M_\Delta$. Since $M_\Delta$
is a multicomplex, it follows that $M_\Delta\subseteq S(\X, (\infty,\a))$.

Lemma III.2.4(a) of \cite{St96} provides a criterion for when a set of
linear forms of $\field[\Delta]$ is an l.s.o.p.
Since $\Delta$ is a pure full-dimensional subcomplex of $\Lambda_d$,
every facet of $\Delta$ is a facet of $\Lambda_d$, and we infer
from this criterion and Lemma \ref{KiKl} that 
$\{g^{-1}(x) : x\in L\}$ is an l.s.o.p.~for $\field[\Delta]$.
As explained in the proof outline, this yields that
$F(M_\Delta)=h(\Delta)$. \endproof

\smallskip\noindent{\it Proof of Lemma \ref{KiKl}: \ }
We need to check that $\det g^{-1}_{H,L}\neq 0$ 
for every subset $H\subset\tilde\X$ corresponding to the
vertex set of a facet of $\Lambda_d$.
For $1\leq i \leq m$, let $H_i=H\cap \X_i$ and $H'_i=H\cap \Y_i$. 
It follows from the definition of $\Lambda_d$ that the set
$H_i\cup H'_i$ corresponds to a face of 
$\Skel_{a_i-1}(\overline{V}_i)$, and hence
$|H_i|+|H'_i|\leq a_i=|\Y_i|$.

Write $I=I_{\X_{m+1}, \X_{m+1}}$ for the identity matrix. 
Since $L\subseteq \X_{m+1}$, we have
$$g^{-1}_{H,L}=\left[\begin{array}{c}
 ABC\\
 C
\end{array} \right]_{H,L} =
\left(\left[\begin{array}{c}
 AB\\
 I
\end{array} \right] C\right)_{H,L},
$$
where $AB$ is a block-diagonal (non-square) matrix whose $i$th block,
$A_iB_i$, has rows indexed by elements of $\X_i$ and columns by elements
of $\Y_i$, $1\leq i\leq m$. By the Cauchy-Binet formula, 
\begin{equation} \label{CauBi}
\det g^{-1}_{H,L}= 
\sum_{\mbox{\tiny $\begin{array}{c}
T\subseteq \X_{m+1},\\ |T|=|H|=d \end{array}$}} 
\det \left[\begin{array}{c}
    AB\\ I \end{array} \right]_{H, T} \cdot
\det C_{T, L}.
\end{equation}

We claim that the determinant of 
{\tiny $\left[\begin{array}{c} AB\\I\end{array} \right]_{H,T}$}
vanishes if 
$T\not\supseteq H'_i$ or if  $|T\cap \Y_i|>|H_i|+|H'_i|$ for some 
$i=1,\ldots, m$. Indeed, in the latter case, the columns of 
{\tiny $\left[\begin{array}{c}AB \\ I\end{array}\right]_{H,T\cap Y_i}$} 
are linearly dependent, while in the former case, $I_{H'_i, T}$
contains a row of zeros. Thus eq. (\ref{CauBi}) yields
\begin{equation}  \label{smaller_sum}
\det g^{-1}_{H,L}=
\sum_{\mbox{\tiny{$\begin{array}{cc}
     T=\cup_{i=1}^m (T_i \cup H'_i)\\
      T_i\subseteq \Y_i-H'_i, \ |T_i|=|H_i| \end{array}$}}} 
\det(A_1B_1)_{H_1,T_1}  \cdot \ldots \cdot \det(A_mB_m)_{H_m,T_m}
\cdot \det C_{T, L}.
\end{equation}
The sum on the right-hand-side of (\ref{smaller_sum}) is not empty:
this is because $|H_i|+|H'_i|\leq a_i=|\Y_i|$ for all $1\leq i \leq m$.
Each of its summands is a non-zero polynomial in the indeterminates 
$(\alpha^j_{xy})_{x,y\in\X_j}$, $1\leq j\leq m+1$, 
and $(\beta^j_{xy})_{x\in\X_j,y\in\Y_j}$, $1\leq j\leq m$.
Moreover, there are no cancellations between different summands 
as for two sets $T\neq T'$, no monomial in 
$(\alpha^{m+1}_{xy})_{x,y\in\X_{m+1}}$ that appears in the expansion of 
$\det C_{T,L}$ appears in the expansion of $\det C_{T', L}$.
Thus $\det g^{-1}_{H,L}\neq 0$. \endproof

To complete the proof of $1.\rightarrow 3.$, it only
remains to verify Lemma \ref{gin}. Its proof is an easy
consequence of the structure of $g$ and
a few known facts about revlex (generic) initial ideals.

\smallskip\noindent{\it Proof of Lemma \ref{gin}: \ }
Fix $1\leq i \leq m$, and let $\I=\I(i)\subset \Field[\tilde\X]$ be
the ideal generated by all squarefree monomials of degree $a_i+1$
supported in $\X_i\cup \Y_i$. Since, as follows from the 
definition of $\Lambda_d$, $\I\subseteq I_{\Lambda_d}$,
and hence $\Init(g\I)\subseteq \Init(gI_{\Lambda_d})$,
to prove the lemma it is enough to show
that $\Init(g\I)$ contains all monomials of degree 
$a_i+1$ supported in $\X_i$. 

By definition of $g$, $g(x)=\sum_{y\in \X_i} (-A_i^{-1})_{xy}y$ if 
$x\in\X_i$, while 
$g(x)=\sum_{y\in \X_i} \beta^i_{xy}y+
\sum_{y\in\X_{m+1}} (C^{-1})_{xy}y$ if $x\in \Y_i$, and so
all monomials
that belong to $\Init(g\I)_{a_i+1}$ are supported in $\X_i\cup\X_{m+1}$. 
Consider a specialization of $g$, 
$\tilde{g}$, obtained by replacing all $(C^{-1})_{xy}$ for 
$y\notin\Y_i$ with zero. Recall that any element of $\X_i$
is $\succ$-larger than any element of $\X_{m+1}\supseteq\Y_i$.
Thus $\tilde{g}$ induces an automorphism of
 $\Field[\X_i\cup \Y_i]$ that is generic
in the sense of Theorem 15.18 of \cite{Eis} and the definition 
following it. Since $\Char \Field=0$, $|\Y_i|=a_i$, and $\I_{a_i+1}$
is the $\Field$-span of all squarefree monomials on $\X_i\cup\Y_i$ 
of degree $a_i+1$, it then follows, e.g., from \cite[Cor.~1.6]{AHH}, 
that $\Init({\tilde g}\I)_{a_i+1}$ is the $\Field$-span of 
 all monomials of degree $a_i+1$ supported in $\X_i$.
Finally, since for every minimal generator of $\I$, and hence also for
every polynomial $\psi\in\I_{a_i+1}$, all 
monomials that appear in the expansion of $g(\psi)-\tilde{g}(\psi)$ 
involve some elements of $\X_{m+1}$, 
and thus are revlex-smaller than any monomial of $\Init({\tilde g}\I)_{a_i+1}$,
we infer that 
$\Init(g\I)_{a_i+1}=\Init({\tilde g}\I)_{a_i+1}$. \endproof  

\begin{question}
Does the implication $1.\to 3.$ continue to hold if 
one relaxes the condition of $h$ being the $h$-vector 
of a $\Q$-CM subcomplex of $\Lambda_d$ to being the 
$h$-vector of a $\field$-CM subcomplex for some field $\field$
of an arbitrary (rather than zero) characteristic?
\end{question}

\section{Compression}

To move from multicomplexes to shellable complexes, we use a generalization of
the compression method of Macaulay \cite{Macaulay}. 
It allows one to replace a general 
multicomplex with a more structured multicomplex having the same $F$-vector. 
This follows \cite{BjFrSt}, whose generalization of Macaulay's theorem to 
colored multicomplexes we extend with our specialized notion of 
``$(0)$-compression,'' and \cite{Br} (which uses a different generalization 
due to Clements and Lindstr\"om, \cite{ClLi}).

As in the introduction, let $\X_0, \X_1, \ldots , \X_m$ be pairwise disjoint 
finite sets of variables and let $\X = \cup_{i=0}^m \X_i$. Let 
$\mathbf{a} = (a_1, a_2, \ldots, a_m)$ be an integer vector and set 
$S = \mult{\X}{\mathbf{a}}{d}$. Fix a total order $\succ$ on $\X$ such that 
$x \succ y$ if $x \in \X_i$ and $y \in \X_j$ with $i <j$. Note that this is 
distinct from the order used in the proof of $1.\rightarrow 3.$ in the 
previous section! In particular, all of the elements of $\X_0$ occur 
\emph{before} those in $\X_i$ for $i>0$, as opposed to after.

Let $M \subseteq S$ be a multicomplex and suppose that $1 \leq i \leq m$. 
Then $M$ is called \emph{$(0, i)$-compressed} if whenever 
$\mu \in M$ and $\mu'$ is a monomial on $\X$ such that 
\begin{enumerate}
\item $\deg(\mu') = \deg(\mu)$, 
\item $\mu_{\X - (\X_0 \cup \X_i)} = \mu'_{\X - (\X_0 \cup \X_i)}$, and
\item  $\mu'_{\X_0 \cup \X_i} \succ \mu_{\X_0 \cup \X_i}$,
\end{enumerate}
then $\mu' \in M$.  If $M$ is $(0, i)$-compressed for each integer $i$ 
with $1 \leq i \leq m$, then $M$ is called \emph{$(0)$-compressed}.

It follows from a result of Mermin and Peeva \cite[Theorem 4.1]{MePe} that if 
$m=1$, so that $S$ is the set of monomials $\mu$ of degree no greater than $d$ 
on $\X_0 \cup \X_1$ with $\deg(\mu_{\X_1}) \leq a_1$, then all possible 
$F$-vectors
of multicomplexes in $S$ are obtained by $(0)$-compressed multicomplexes. 
In this case, it simply means that for $M$ a multicomplex in $S$, the set 
containing the first $F_i(M)$ elements of $S$ in degree $i$ for each $i$ 
is a multicomplex. We may use this result to obtain the following theorem, 
which specializes to $3.\rightarrow 4.$

\begin{theorem} \label{compressiontheorem} 
Let $M \subseteq \mult{\X}{\mathbf{a}}{d}$ be a multicomplex. Then there exists 
a multicomplex $M' \subseteq \mult{\X}{\mathbf{a}}{d}$ such that $M'$ is 
$(0)$-compressed and $F(M') = F(M)$.
\end{theorem}
\proof
For $M$ a multicomplex in $S$ and $1 \leq i \leq m$,  
let $M(i) = \{\mu_{\X - (\X_0 \cup \X_i)}  : \mu \in M\}$. For $\nu \in M(i)$, 
let $M_{\nu}$ be the set of monomials $\mu$ in $\X_0 \cup \X_i$ such that 
$\mu\nu \in M$. Observe that $M_{\nu}$ is a multicomplex. In particular, 
by \cite{MePe} the set $M'_{\nu}$ containing the first $F_j(M_{\nu})$ elements 
of $\mult{\X_0 \cup \X_i}{a_i}{d}$ of degree $j$ for each $j\leq d$ is a 
multicomplex.

Define $C_i(M) := \cup_{\nu \in M(i)}M'_{\nu}\nu$. Then 
$C_i(M)$ is $(0, i)$-compressed and $F(C_i(M)) = F(M)$.
We next show that $C_i(M)$ is a multicomplex.

Suppose $\mu\nu \in C_i(M)$ where $\nu$ is supported in $\X - (\X_0 \cup \X_i)$ 
and $\mu$ is supported in $\X_0 \cup \X_i$. Any divisor of $\mu\nu$ is of the 
form $\mu'\nu'$ where $\nu'$ and $\mu'$ are supported in 
$\X - (\X_0 \cup \X_i)$ and $\X_0 \cup \X_i$, respectively.
Then $\nu' | \nu$, and as $M$ is a multicomplex, $M_{\nu} \subseteq M_{\nu'}$. 
Thus $M'_{\nu} \subseteq M'_{\nu'}$. In particular, since $\mu \in M'_{\nu}$, 
$\mu \in M'_{\nu'}$. Finally, since $\mu' | \mu$ and $M'_{\nu'}$ is a multicomplex, 
$\mu' \in M'_{\nu'}$ and $\mu'\nu' \in C_i(M)$.

Let $\minv_{M} = \prod_{\mu \in M} \mu$. It is immediate from the definition 
that $\deg(\minv_{M}) = \deg(\minv_{C_i(M)})$ and 
$\minv_{C_i(M)} \succeq \minv_M$, with equality if and only if $C_i(M) = M$, 
in which case $M$ is $(0,i)$-compressed.

To complete the proof we apply $C_i$ repeatedly to obtain our $(0)$-compressed 
multicomplex. Let $M = M_0$ and inductively define $M_{i +1} = C_{j}(M_i)$, 
where $1 \leq j \leq m$ and $i+1 \equiv j \mod m$. Then 
$\minv_{M_0} \preceq \minv_{M_1} \preceq \minv_{M_2} \preceq \cdots$. As this 
sequence cannot increase indefinitely there must be some $k$ such that 
$\minv_{M_k} = \minv_{M_j}$ for all $j > k$. In particular $C_i(M_k) = M_k$ 
for each $i = 1, \ldots, m$, and so $M_k$ is a $(0)$-compressed multicomplex.  
\endproof

\section{From multicomplexes to shellable complexes}
We are now in a position to lay the groundwork for our proof of the implication 
$4.\rightarrow 2.$ in Theorem \ref{main}. Throughout this section and the next 
fix  $\V = \V_1 \cup \V_2 \cup \cdots \cup \V_m$, 
$\mathbf{a} = (a_1, a_2, \ldots, a_m)$ a non-negative integer vector with 
$a_i \leq |\V_i|$, and $1 \leq d \leq \sum a_i$. Let $n = |\V|$. For brevity let 
$\Lambda = \comp{\V_1, \V_2, \ldots \V_m}{\mathbf{a}}{d}$. Let $\fac(\Lambda)$ 
denote the set of facets of $\Lambda$. Order the elements of $\V$ such that if 
$v \in \V_i$ and $v' \in \V_j$ with $i < j$, then $v \succ v'$. 
Let  $\X = \X_0 \cup \X_1 \cup \cdots \cup \X_m$, where 
$|\X_i| = |\V_i| - a_i$ for $i = 1, \ldots, m$ and 
$|\X_0| = (\sum_{i=1}^m a_i) -d$. Order the elements of $\X$ such that 
$x \succ x'$ if $x \in \X_i$, $x' \in \X_j$ and $i < j$, and
let $S = \mult{\X}{\mathbf{a}}{d}$.

\medskip\noindent{\bf Outline of $4.\rightarrow 2.$}
The structure of the proof is as follows: first, we let $\Lex$ be the listing of 
the facets of $\Lambda$ in the revlex order, i.e, 
$\Lex= (\tau_1, \tau_2, \ldots)$ where $\tau_1 \succ \tau_2 \succ \cdots$, and 
show in Lemma \ref{shell}
that this is a shelling of $\Lambda$. In doing so we explicitly determine 
its restriction function $\Res_{\Lex}$. We then establish the following theorem, 
which we prove in Section \ref{bijectionproof}. The implication 
$4.\rightarrow 2.$ in Theorem \ref{main} follows (see Corollary \ref{4-->2}).

\begin{theorem} \label{bijection} There exists a bijection 
$\Phi : \fac(\Lambda) \rightarrow S$ such that for each $\tau \in \fac(\Lambda)$,
 $\deg(\Phi(\tau)) = |\Res_{\Lex}(\tau)|$, and whenever 
$\Res_{\Lex}(\tau) \nsubseteq \gamma \subset \tau$, there is a facet 
$\tau'$ of $\Lambda$, divisor $\mu''$ of $\mu = \Phi(\tau)$ and integer 
$i$ such that 
\begin{enumerate}[(a)]
\item $\gamma \subset \tau'$,
\item $\tau' \succ \tau$, 
\item $\deg(\Phi(\tau')) = \deg(\mu'')$, 
\item $\Phi(\tau')_{\X_0 \cup \X_i} \succeq \mu''_{\X_0 \cup \X_i} $,
\item $\Phi(\tau')_{\X -  (\X_0 \cup \X_i)} = \mu''_{\X -  (\X_0 \cup \X_i)}$.
\end{enumerate}
\end{theorem}

\begin{corollary} \label{4-->2}
Let $M$ be a $(0)$-compressed multicomplex in $S$ and 
$\Gamma$  the simplicial complex whose facets are exactly 
$\Phi^{-1}(M)$. Then $\Gamma$ is shellable, and $h(\Gamma) = F(M)$.
\end{corollary}

\proof Since  $\deg(\Phi(\tau)) = |\Res_{\Lex}(\tau)|$, 
by eq.~(\ref{h-shellable}) it suffices to show that by putting the elements of 
$\Phi^{-1}(M)$ in revlex order we obtain a shelling of $\Gamma$ whose restriction
 function is simply the restriction of $\Res_{\Lex}$ to $\Phi^{-1}(M)$.

First, we note that there is no facet $\tau' \succ \tau$ of $\Gamma$  
containing $\Res_{\Lex}(\tau)$: indeed, by Lemma~\ref{shell},
$\Lex$  is a shelling of $\Lambda$ and each facet of $\Gamma$ is a facet of $\Lambda$. 
On the other hand, suppose $\Res_{\Lex}(\tau) \nsubseteq \gamma \subset \tau$ 
for some facet $\tau$ of $\Gamma$. Let $\tau'$ and $\mu''$ be the facet of 
$\Lambda$ and divisor of $\mu = \Phi(\tau)$ given by Theorem \ref{bijection}. 
Since $\mu''$ divides $\mu$, and $\mu \in M$, $\mu'' \in M$. Then, by 
properties $(c)$, $(d)$ and $(e)$, $\Phi(\tau')$ must be in $M$, as $M$ is 
$(0)$-compressed. In particular, $\tau'$ is a facet of $\Gamma$, 
contains $\gamma$ and occurs earlier than $\tau$ in the revlex order. 
Thus $\Res_{\Lex}(\tau)$ is the unique minimal face of 
$\overline{\tau} - 
(\cup_{\tau' \succ \tau, \,  \tau' \in \Phi^{-1}(M)}\overline{\tau'})$.
\endproof

\medskip\noindent{\bf The Shelling.}
 Recall our definition of the restriction function, 
$\Res_{\Lex}(\tau) = \{v\in\tau \ : \ 
\tau-v \subseteq \tau' \mbox{ for some } \tau' \succ \tau \}$. 
For $v \in \tau$, any other facet of $\Lambda$ containing $\tau - v$ is of the form 
$(\tau - v) \cup w$ for some $w \notin \tau$, and occurs earlier than $\tau$ if and 
only if $w \succ v$. In other words, the elements of 
$\Res_{\Lex}(\tau)$ are precisely those $v \in \tau$ such that we may ``swap'' $v$ 
for some $w \succ v$, $w \notin \tau$,  without leaving $\Lambda$.

Which vertices are these? To start, there must be  some $w \notin \tau$ 
with $w \succ v$. Now, $w \in \V_i$ for some $i$; if $| \tau \cap \V_i| = a_i$, 
$(\tau - v) \cup w$ may contain too many elements of $\V_i$ and thus not be in 
$\Lambda$. To distinguish such $\V_i$, we define 
$\FL(\tau) = \{ i : |\tau \cap \V_i| = a_i \}$; 
that is, $\FL(\tau)$ is the collection of indices of sets $\V_i$ from which $\tau$ 
contains the maximum allowed number of vertices 
(so $\tau$'s intersection with these sets is ``full''). Let $\Gap(\tau)$ be the 
earliest element of $\V$ which is in neither $\tau$ nor any $\V_i$ with 
$i \in \FL(\tau)$, if such an element exists; then if 
$\Gap(\tau) \succ v$, we may take $w = \Gap(\tau)$.  
The set of such $v$'s forms a ``tail'' of $\tau$, 
$$\tail(\tau) : =  \{ v \in \tau : v  \prec \Gap(\tau) \}.$$
If every element of $\V - \tau$ is contained in some $\V_i$ with $i \in \FL(\tau)$, 
set $\tail(\tau) = \emptyset$.

On the other hand, suppose the only $w \succ v$ which are not in $\tau$ are in 
$\V_i$ where $i \in \FL(\tau)$. Then $(\tau - v) \cup w$ is in $\Lambda$ 
if and only if $w$ and $v$ are both in $\V_i$. This case occurs when $v$ 
occurs after the first element of $\V_i$ not in $\tau$. Hence we let 
$\fgap(\tau, i) = \max_{\succ} (\V_i - \tau )$ be the ``first gap'' in 
$\V_i \cap \tau$, for each $i \in \FL(\tau)$ such that  $\V_i \cap \tau \neq \V_i$. Let 
$$\up(\tau) := \{ v \in \V : 
v \in \V_i,\  i \in \FL(\tau), \ 
\V_i \cap \tau \neq \V_i, \text{ and } v \prec \fgap(\tau, i) \}$$ 
be the ``lower'' part of $\tau$ in the segments of the partition it intersects 
maximally. Then $\Res_{\Lex}(\tau) = \up(\tau) \cup \tail(\tau)$.

\begin{example} \label{R-example}
{\rm Suppose $\V = \V_1 \cup \V_2 \cup \V_3$, where 
$\V_1 = \{ v^1_1, v^1_2, v^1_3 , v^1_4 \}$, 
$\V_2 = \{ v^2_1, v^2_2, v^2_3 \}$,  and $\V_3 = \{ v^3_1,v^3_2, v_3^3 \}$, 
and let $\a = (2,2,1)$.  Order $\V$ such that $v^s_j \succ v^t_i$ if either 
$s <t$ or $s=t$ and $j < i$. We depict $\V$ pictorially, with columns corresponding 
to the parts $\V_i$ of our partition. The order reads top to bottom, left to right.
\begin{center}
\includegraphics{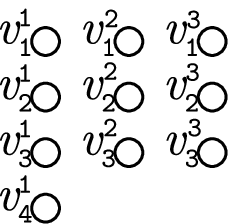}
\end{center}
Taking $d = 4$, $\tau = \{v^1_1, v^1_4 , v^2_2, v^3_2 \}$ is a facet of $\Lambda$.
\begin{center}
\includegraphics{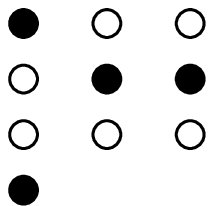}
\end{center}
The first and third columns are ``full''; i.e., $\FL(\tau)  =  \{1 , 3 \}$. 
So $\Gap(\tau)$ is the first missing element of $\V_2$, $v^2_1$, and 
$\tail(\tau)  =  \{ v^2_2, v^3_2 \}$.
\begin{center}
\includegraphics{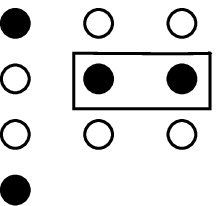}
\end{center}
On the other hand, $\up(\tau)  =  \{ v^1_4, v^3_2 \}$, the lower elements of the 
full columns.
\begin{center}
\includegraphics{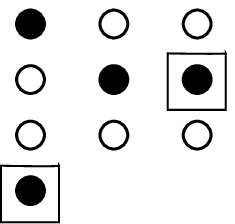}
\end{center}
Now consider $\up(\tau) \cup \tail(\tau)$:
\begin{center}
\includegraphics{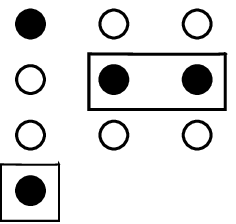}
\end{center}
Observe that this is the unique minimal subset of $\tau$ that is 
contained in no facet occuring earlier in the reverse lexicographic order.}
\end{example}

\begin{lemma}  \label{shell}
For every facet $\tau$ of $\Lambda$,  
$\Res_{\Lex}(\tau)$ is the unique 
minimal face of $\overline{\tau} - (\cup_{\tau' \succ \tau}\overline{\tau'})$. 
In particular, $\Lex$ is a shelling of $\Lambda$.
\end{lemma}
\proof
Let $\tau$ be a facet of $\Lambda$. As $\Res_{\Lex}(\tau) = \{v\in\tau \ : \ 
\tau-v \subseteq \tau' \mbox{ for some } \tau' \succ \tau \}$, 
any face of $\tau$ not containing $\Res_{\Lex}(\tau)$ is in 
$(\cup_{\tau' \succ \tau}\overline{\tau'})$, so every face
of $\overline{\tau} - (\cup_{\tau' \succ \tau}\overline{\tau'})$ 
contains $\Res_{\Lex}$. Thus it suffices to show that 
$\Res_{\Lex} \in \overline{\tau} - (\cup_{\tau' \succ \tau}\overline{\tau'})$; 
that is, there is no facet $\tau' \succ \tau$ of $\Lambda$ 
containing $\Res_{\Lex} = \up(\tau) \cup \tail(\tau)$.

Suppose, in order to obtain a contradiction, that there is such a facet $\tau'$. 
As $\tau' \succ \tau$, there is vertex $v$ such that $v \in \tau$,
$v \notin \tau'$, and $\tau$ and $\tau'$ agree on all vertices after $v$. 
Then, since $\tau$ and $\tau'$ contain the same number of elements, 
there must be some $w \succ v$ such that $w \in \tau'$ and $w \notin \tau$.

Can it happen that $v \prec \Gap(\tau)$? No, since then we would have
$v \in \tail(\tau) \subseteq \tau'$, 
a contradiction. Thus $v \succ \Gap(\tau)$. Now let $L$ be the set of elements 
in $\V$ strictly less than $v$, and $U$ be the set of elements of $\V$ 
strictly greater than $v$. In particular, we see that 
$| \tau' \cap U| = | \tau \cap U| + 1$. 
Hence there must be some $i$ such that 
$|\V_i \cap \tau \cap U| < |\V_i \cap \tau' \cap U|$.

For such an $i$, $\tau$ cannot contain all of $\V_i \cap U$, and so since 
$v \succ \Gap(\tau)$, we must have $i \in \FL(\tau)$. 
If $v \in \V_i$, then there is an element of $\V_i$ that
is not in $\tau$ and is greater than $v$, and hence as $i \in \FL(\tau)$, 
$v \in \up(\tau)$. But then $v \in \tau'$, a contradiction. 
Thus $v \notin \V_i$. In particular, 
$\tau \cap \V_i = (\tau \cap \V_i \cap L) \cup (\tau \cap \V_i \cap U)$,
and we obtain that
$$
|\tau' \cap \V_i|  =  |\tau' \cap \V_i \cap L| + |\tau' \cap \V_i \cap U|                   
                   >  |\tau \cap \V_i \cap L| + |\tau \cap \V_i \cap U|              
                          =  a_i,     
$$
which is again a contradiction. Therefore, no such $\tau'$ may exist.
\endproof

\section{The Bijection} \label{bijectionproof}
In this section we complete the proof of our main 
theorem --- Theorem \ref{main}. To do this, it only remains to
verify Theorem \ref{bijection} on the existence of $\Phi$.
Our bijection $\Phi$ is a generalization of one used in \cite{BjFrSt},
and is similarly built starting from a map corresponding to what in our
notation is the case $m =1$, $d=a_1$ (i.e., $\V = \V_1$ and $\X = \X_1$).
Let $V$ be a finite set of vertices, $a \leq |V|$ a positive integer and $X$
a set of $|V| - a$ variables. Put some total order $\succ$ on $V$ and label
its elements $v_1\succ v_2 \succ \cdots$ accordingly. For $\tau$ an $a$-subset
of $V$, we may write
$\tau = \{ v_1, v_2, \ldots, v_t, v_{i_1}, v_{i_2}, \ldots, v_{i_s} \}$
where $t,s\geq 0$, $i_1 > t+1$, and $t + s = a$. Then define
\begin{equation*}                                                          
\phi(V, X)(\tau) = x_{i_1 -(t+1)}x_{i_2 - (t+2)}\cdots x_{i_s - (t+s)}.             
\end{equation*}
Thus $\phi(V,X)$ maps the set of $a$-subsets of $V$ into the set of monomials
in $X$ with degree no greater than $a$. On the other hand, for $\mu$ a monomial
in $X$ with degree no greater than $a$, we may write $\mu = x_{i_1}\cdots x_{i_s}$
with $i_1 \leq i_2 \leq \ldots \leq i_s$ and define
\begin{eqnarray*} 
\sigma(X,V)(\mu) & = & \{v_1, v_2, \ldots, v_{a-s} \},\\ 
\rho(X,V)(\mu) & = & \{v_{i_1 + a-(s-1)}, v_{i_2 +a-(s-2)}, \ldots , 
  v_{i_s + a} \},\\     
\psi(X,V)(\mu) & = & \sigma(X,V)(\mu) \cup \rho(X,V)(\mu).
\end{eqnarray*}
Then $\psi(X,V)$ maps the set of monomials in $X$ of degree no greater than $a$
into the set of $a$-subsets of $V$. It is easy to see that $\phi(V,X)$ and 
$\psi(X,V)$
are inverse to each other. The usefulness of these bijections is explained by 
the following lemma.

\begin{lemma} \label{littlephi} Let $\tau, \tau'$ be $a$-subsets of $V$
taken to $\mu, \mu'$, respectively, by $\phi(V, X)$.
If $\tau' \succ \tau$ and $\deg(\mu') \leq \deg(\mu)$, then there is a divisor
$\mu''$ of $\mu$ such that $\deg(\mu'') = \deg(\mu')$ and $\mu' \succeq \mu''$.
\end{lemma}
\proof Write
$$ \mu =  x_{i_1}\cdots x_{i_s} \quad \mbox{and} \quad                                 
\mu' =  x_{i'_1}\cdots x_{i'_{s'}},$$
with $i_1 \leq i_2 \leq \ldots \leq i_s$ and
$i'_1 \leq i'_2 \leq \ldots \leq i'_{s'}$,
so that
\begin{eqnarray*}    
\tau & = & \{v_1, v_2, \ldots, v_{a-s} \}  
\cup \{v_{i_1 + a-(s-1)}, v_{i_2 +a-(s-2)}, \ldots , v_{i_s + a} \},\\              
\tau' & = & \{v_1, v_2, \ldots, v_{a-s'} \}                                         
 \cup \{v_{i'_1 + a-(s'-1)}, v_{i'_2 +a-(s'-2)}, \ldots , v_{i'_{s'} + a}\}.      
\end{eqnarray*}
Now, as $\tau' \succ \tau$, there is some $v_j$ such that
$v_j \in \tau$, $v_j \notin \tau'$, and $\tau$, $\tau'$ agree on all vertices
after $v_j$. In particular, $j > a-s' \geq a-s$ (as $v_j \notin \tau'$
and $s' \leq s$). Thus $j = i_{s-r} + (a-r)$, for some $0\leq r<s$.
Since $\tau$ and $\tau'$ agree after $v_j$, the last $r$
elements of each are the same, i.e,  $r\leq s'$ and for $l = 1, \ldots r$,
\begin{eqnarray*}                                                                   
i_{s-r+l} + a-(r-l) &=& i'_{s'-r + l} + a-(r-l), \text{ and so}\\    
i_{s-r+l} & = &i'_{s'-r+l}.                                                         
\end{eqnarray*}
 On the other hand, as $v_j \notin \tau'$, in the case of $r<s'$, we also have 
 $$i'_{s'-r} + a-r < j=  i_{s-r}+a-r, \quad \mbox{hence} \quad i'_{s'-r}< i_{s-r}.$$    
It follows that if we let $\mu''= \prod_{l = s-s'+1}^{s}x_{i_l}$,
then $\deg(\mu'') = s' = \deg(\mu')$ and $\mu' \succeq \mu''$.
\endproof

We now build $\Phi$ from $\phi$. Let $\tau \in \fac(\Lambda)$, $1 \leq i \leq m$, 
and consider $\tau \cap \V_i$. We would like to apply $\phi(\V_i, \X_i)$, 
but we note that $\tau \cap \V_i$ may contain fewer than $a_i$ elements, 
in which case it is not in the domain. So define  $\fll{\tau}{i}{a_i}$ 
to be the reverse lexicographically first $a_i$-subset of $\V_i$ containing 
$\tau \cap \V_i$, and set $\Phi_i(\tau):=\phi(\V_i, \X_i)(\fll{\tau}{i}{a_i})$.

Next, let $\V[\tau]$ be the subset of $\V$ which contains the first 
$a_i - \deg(\Phi_i(\tau))$ elements of $\V_i$ for $1 \leq i \leq m$, 
and let $\tau[0] = \tau \cap \V[\tau]$. Our aim is to apply $\phi(\V[\tau], \X_0)$
to  $\tau[0]$. To do so we must show that $\tau[0]$ is in the domain. 
This is done in Lemma \ref{domain}. We then define 
 $\Phi_0(\tau) := \phi(\V[\tau], \X_0)(\tau[0])$, and finally set
\begin{equation*}    
\Phi(\tau) := \prod_{i=0}^m \Phi_i(\tau). 
\end{equation*}

\begin{example}{\rm Return to the $\Lambda$ and $\tau$ of Example \ref{R-example}. 
The corresponding set of variables is $\X = \X_0 \cup \X_1 \cup \X_2 \cup \X_3$ 
where $|\X_1| = 2$, $|\X_2| = 1$, $|\X_3| = 2$ and $|\X_0| = 1$. 
Label these $\X_0 = \{ w\}$, $\X_1 = \{ x_1, x_2\}$, $\X_2 = \{ y\}$, 
$\X_3 = \{ z_1, z_2\}$, ordered accordingly.                 
Then $\fll{\tau}{1}{2} = \tau \cap \V_1 = \{v^1_1, v^1_4 \}$, and so 
$\Phi_1(\tau) = x_2$. Similary 
$\fll{\tau}{3}{1} = \tau \cap \V_3 = \{v^3_2 \}$, and hence
$\Phi_3(\tau) = z_1$. On the other hand, 
$\fll{\tau}{2}{2} = \{v^2_1, v^2_2 \}$, so $\Phi_2(\tau) = 1$.               
Thus $\V[\tau] = \{v^1_1, v^2_1, v^2_2 \}$, 
so that $\tau[0]= \{v^1_1, v^2_2 \}$   
and $\Phi_0(\tau) = w$. Putting  all this together, we obtain that $\Phi(\tau) = wx_2z_1$. 
Note that $\deg(\Phi(\tau)) = 3 = |\Res_{\Lex}(\tau)|$.}
\end{example}

To complete our definition of $\Phi$, it only remains to check the followng.
\begin{lemma} \label{domain}
$|\V[\tau]| - |\tau[0]| = |\X_0|$.
\end{lemma}
\proof 
It follows from the definition of $\Phi_i$ that $\fll{\tau}{i}{a_i}$ consists of
its initial segment and $\deg(\Phi_i(\tau))$ elements outside of it. These last 
$\deg(\Phi_i(\tau))$ elements all belong to $\tau$, while the initial segment, 
whose length is $a_i-\deg(\Phi_i(\tau))$, contains $|\tau[0]\cap \V_i|$ elements 
of $\tau$ and (possibly) a few added elements. Thus
$$d   =  |\tau|
    =  |\tau[0]| + \sum_{i=1}^m \deg(\Phi_i(\tau)),   \quad \mbox{and so}$$
$$\left|\V[\tau]\right| = \sum_{i=1}^m \left(a_i - \deg \Phi_i(\tau)\right) 
           =  (\sum_{i=1}^m a_i) -\left( d-|\tau[0]| \right)
    = |\X_0| + |\tau[0]|,
$$ 
as required.   \endproof

To complete the proof of our main theorem, it only remains to show
that $\Phi$ satisfies the conditions of Theorem \ref{bijection}. 
We start by showing that $\Phi$ is a bijection. To do so,
 we explicity construct its inverse.

Let $\mu \in S$. Then for $1 \leq i \leq m$, $\deg(\mu_{\X_i}) \leq a_i$, 
so we may define $\Psi_i(\mu) = \rho(\X_i, \V_i)(\mu_{\X_i})$. 
Next, let $\V[\mu]$ be the subset of $\V$ containing the first 
$a_i - \deg(\mu_{\X_i})$ elements of $\V_i$ for each $i = 1, \ldots ,m$.  

Notice that $\mu_{\X_0}$ has degree no greater than 
$d - \sum_{i =1}^m \deg(\mu_{\X_i})$, while $\psi(\X_0, \V[\mu])$ 
takes  monomials on $\X_0$ of degree no greater than 
$$
|\V[\mu]| - |\X_0| =  
\sum_{i=1}^m( a_i - \deg(\mu_{\X_i})) - \left((\sum_{i=1}^m a_i) - d\right)\\
=  d - \sum_{i =1}^m \deg(\mu_{\X_i})
$$
to subsets of $\V[\mu]$ of size $d - \sum_{i =1}^m \deg(\mu_{\X_i})$. 
Thus if we let $\Psi_0(\mu) = \psi(\X_0, \V[\mu])(\mu_{\X_0})$, and set
 \begin{equation*}
\Psi(\mu): = \cup_{i=0}^m \Psi_i(\mu),
\end{equation*}
then $\Psi(\mu)$ is a $d$-subset of $\V$. Next note that for $1 \leq i \leq m$,
\begin{eqnarray*}
|\Psi(\mu) \cap \V_i| & = & |\Psi(\mu) \cap \V_i \cap \V[\mu]| + 
|\Psi(\mu) \cap \V_i \cap (\V - \V[\mu])|\\
& = &  |\Psi(\mu) \cap \V_i \cap \V[\mu]| + |\Psi_i(\mu)|\\
& \leq & |\V_i \cap \V[\mu]| + |\Psi_i(\mu)|\\
& = & a_i - \deg(\mu_{\X_i}) + \deg(\mu_{X_i}) = a_i.
\end{eqnarray*}
In other words, $\Psi(\mu)$ is a facet of  $\Lambda$.

\begin{example} {\rm Return to the $\Lambda$ and $S$ of the previous two examples. 
Let $\mu = x_2z_1$. Then $\Psi_1(\mu) = \{v^1_4\}$, $\Psi_2(\mu) = \emptyset$ 
and $\Psi_3(\mu)= \{v^3_2\}$. Furthermore, $\V[\mu] = \{v^1_1, v^2_1, v^2_2 \}$, 
and $\Psi_0(\mu) = \{ v^1_1, v^2_1 \}$. Thus $\Psi(\mu) = \{v^1_1, v^1_4, v^2_1, v^3_2 \}$. 
For comparison, if we instead take $\nu = wx_2z_1$ we obtain the same 
$\Psi_1, \Psi_2, \Psi_3$ and $\V[\nu]$, but $\Psi_0(\nu) =  
\{ v^1_1, v^2_2 \}$, so $\Psi(\nu) =  \{v^1_1, v^1_4 , v^2_2, v^3_2 \}$, 
that is, $\tau$ from the earlier examples (as we should hope if $\Psi$ is to invert $\Phi$).}
\end{example}

\begin{lemma} \label{inverse} 
$\Psi$ and $\Phi$ are inverse.
\end{lemma}
\proof
It follows from the implication $1.\to 3.$ of Theorem \ref{main} (proved in Section 3), 
that the $h$-vector of $\Lambda$ is the $F$-vector of a sub-multicomplex of $S$. Thus
$$
|\fac(\Lambda)| =f_{d-1}(\Lambda)=\sum_{i=0}^d h_i(\Lambda) \leq \sum_{i=0}^d F_i(S) = |S|.
$$
Since $\Phi: \fac(\Lambda)\to S$ and $\Psi: S\to\fac(\Lambda)$, to prove the lemma
it suffices to show that the composition $\Phi \circ \Psi : S\to S$ is the identity map. 

Let $\mu \in S$ and let $\tau = \Psi(\mu)$. Observe that for $1 \leq i \leq m$, 
the earliest $a_i$-subset of $\V_i$ containing $\tau \cap \V_i$  is that which 
contains the $\deg(\mu_{\X_i})$ elements of $\Psi_i(\mu)$ and all 
$a_i - \deg(\mu_{\X_i})$ elements of $\V[\mu] \cap \V_i$. Hence
$$\fll{\tau}{i}{a_i} = \Psi_i(\mu) \cup \sigma(\X_i, \V_i)(\mu_{\X_i})
 =  \rho(\X_i, \V_i)(\mu_{\X_i}) \cup \sigma(\X_i, \V_i)(\mu_{\X_i})
 =  \psi(\X_i, \V_i)(\mu_{\X_i}), 
$$
so that 
\begin{equation}    \label{i}
\Phi_i(\tau) = \phi(\V_i, \X_i)(\fll{\tau}{i}{a_i})
 =   \phi(\V_i, \X_i)(\psi(\X_i, \V_i)(\mu_{\X_i}) ) 
 =  \mu_{\X_i},
\end{equation}
as $\phi(\V_i, \X_i)$ and $\psi(\X_i, \V_i)$ are inverses. 
Futhermore, this implies that $\V[\tau] = \V[\mu]$. Thus
$$\tau \cap \V[\tau]  =   \tau \cap \V[\mu]
=\Psi(\mu)  \cap \V[\mu] =  \Psi_0(\mu)
 =  \psi(\X_0, \V[\mu])(\mu_{\X_0}),
$$
and so 
\begin{equation}  \label{0}
\Phi_0(\tau)  =  \phi(\V[\tau], \X_0)(\tau \cap \V[\tau])
  =  \phi(\V[\mu], \X_0)(\psi(\X_0, \V[\mu])(\mu_{\X_0}))
  =  \mu_{\X_0}.
\end{equation}
Equations (\ref{i}) and (\ref{0}) imply that $\Phi(\Psi(\mu))=\Phi(\tau) = \mu$,
and the assertion follows.  
\endproof

\begin{lemma} For each facet $\tau$ of $\Lambda$, $\deg(\Phi(\tau)) = |\Res_{\Lex}(\tau)|$.
\end{lemma}
\proof
Let $\tau$ be a facet of $\Lambda$ and let $\mu = \Phi(\tau)$ (so $\tau = \Psi(\mu)$). 
Observe that 
$$\deg(\mu) = \sum_{i=0}^m \deg(\mu_{\X_i})
 =  |\rho(\X_0, \V[\mu])(\mu_{\X_0})| + \sum_{i=1}^m |\rho(\X_i, \V_i)(\mu_{\X_i})|.
$$
Thus it will suffice to show that 
$\Res_{\Lex}(\tau) = 
\rho(\X_0, \V[\mu])(\mu_{\X_0}) \cup (\cup_{i=1}^m \rho(\X_i, \V_i)(\mu_{\X_i})).$

In the case of $v \in \rho(\X_i, \V_i)(\mu_{\X_i})$ for some $i$, $1 \leq i \leq m$, 
it follows from the definition of $\sigma$ that the element of $\V_i$ in position 
$(a_i - \deg(\mu_{X_i}) + 1)$ is not in $\sigma(\X_i, \V_i)(\mu_{\X_i})$. 
This element is also not in $\V[\mu]$. In particular, this is an element of $\V_i$ 
(the set containing $v$) that occurs before $v$ and is not in $\tau$. 
Hence either $v \in \up(\tau)$ or $v \in \tail(\tau)$, and so $v\in \Res_{\Lex}(\tau)$.

In the case of $v \in \rho(\X_0, \V[\mu])(\mu_{\X_0})$, there is an element 
$w \in \V[\mu]$ that is not in $\tau$ and occurs before $v$. Suppose, in order 
to obtain a contradiction, that $w \in \V_i$ where $i \in \FL(\tau)$. Then 
$\fll{\tau}{i}{a_i} = \tau \cap \V_i$, so $\tau$ contains the first 
$a_i - \deg(\mu_{\X_i})$ elements of $\V_i$. But these are all of the elements of 
$\V[\tau] \cap \V_i$, including $w$, a contradiction. 
Hence $v \in \tail(\tau)\subseteq \Res_{\Lex}(\tau)$. 

In the case of $v \in \tau$ and 
$v \notin \rho(\X_0, \V[\mu])(\mu_{\X_0}) \cup 
(\cup_{i=1}^m \rho(\X_i, \V_i)(\mu_{\X_i}))$, we have
$v \in \sigma(\X_0, \V[\tau])(\mu_{\X_0})$; in other words, $v$ is in the 
initial segment of $\tau$ in $\V[\tau]$.
In particular, $v \in \V_j$ for some $j$, and $\tau$ contains every element of $\V_j$ 
occurring before $v$, yielding that $v \notin \up(\tau)$. 
Consider $i < j$. Then $\tau$ contains every element of $\V[\tau] \cap \V_i$. 
As $\tau$ contains $\deg(\mu_{\X_i})$ elements 
of $\V_i$ outside of $\V[\tau]$, $\tau$ in total contains 
$a_i - \deg(\mu_{\X_i}) + \deg(\mu_{\X_i}) = a_i$ elements of $\V_i$, and so 
$i \in \FL(\tau)$. Therefore every  vertex occurring before $v$ is either in $\V_i$ 
with $i \in \FL(\tau)$ or in $\V_j$, and hence in $\tau$. Thus $v \notin \tail(\tau)$, and
we infer that $v \notin \Res_{\Lex}(\tau)$. 
\endproof

We next describe the facet $\tau'$ in Theorem \ref{bijection}. For $\tau$ a facet of 
$\Lambda$ and $\Res_{\Lex}(\tau) \nsubseteq \gamma \subset \tau$, there is an element 
$v$  of $\Res_{\Lex}$ which is not in $\gamma$. If $v \in \tail(\tau)$, 
let $w = \Gap(\tau)$. Otherwise, $v \in \up(\tau)$; 
let $w = \fgap(\tau, i)$, where $i$ is the index of the set $\V_i$ containing $v$. 
Then $\bk{\tau}{\gamma}:= (\tau - v) \cup w$ is a facet of 
$\Lambda$ containing $\gamma$ and $\bk{\tau}{\gamma}\succ \tau$. 
Let $\mu = \Phi(\tau)$ and $\mu' = \Phi(\bk{\tau}{\gamma})$.

\begin{lemma} \label{R-back} 
For $\tau$ and $\gamma$ as above, 
$|\Res_{\Lex}(\bk{\tau}{\gamma})| \leq |\Res_{\Lex}(\tau)|$. 
In particular, $\deg(\mu') \leq \deg(\mu)$.
\end{lemma}
\proof Let $v$ and $w$ be as in the definition of $\bk{\tau}{\gamma}$, 
and suppose $v \in \V_i$. Since $v \in \Res_{\Lex}(\tau)$ but 
$v \notin \Res_{\Lex}(\bk{\tau}{\gamma})$, it will suffice to show that 
$\Res_{\Lex}(\bk{\tau}{\gamma}) \subset \Res_{\Lex}(\tau) \cup w$. 
Let $u \in \Res_{\Lex}(\bk{\tau}{\gamma})$, $u \neq w$ (so, in particular, 
$u \in \tau$). There are two possible cases.

\smallskip\noindent{\bf Case 1:} $u \in \tail(\bk{\tau}{\gamma})$. 
If also $u \in \tail(\tau)$, then $u \in \Res_{\Lex}(\tau)$, and we are done.
Thus assume without loss of generality that $u \notin \tail(\tau)$. 
Then either $\Gap(\bk{\tau}{\gamma}) \in \tau$ or $\Gap(\bk{\tau}{\gamma}) \in \V_j$ 
for some $j \in \FL(\tau)$.

In the case of $\Gap(\bk{\tau}{\gamma}) \in \tau$, we have $\Gap(\bk{\tau}{\gamma}) = v$. 
(This follows from the observation that  $v$ is the only element of $\tau$ 
that is not in $\bk{\tau}{\gamma}$). 
Thus $w \succ v \succ u$. Now, $w \in \V_j$ for some $j$; 
notice that $j\in\FL(\tau)$ or otherwise $u$ would be in $\tail(\tau)$.  
Then by our definition of $w$, $j = i$. 
Since $\bk{\tau}{\gamma}$ must then contain the same number of elements in $\V_i$
as $\tau$, we obtain that $i \in \FL(\bk{\tau}{\gamma})$, a contradiction to
the fact that $v = \Gap(\bk{\tau}{\gamma})$.

In the case of  $\Gap(\bk{\tau}{\gamma}) \in \V_j$ for some $j \in \FL(\tau)$,
we have, by definition of  $\Gap(\bk{\tau}{\gamma})$, 
that $j \notin \FL(\bk{\tau}{\gamma})$. 
Thus $v \in \V_j$ and $w \notin \V_j$ (as this is the only way in which 
$\bk{\tau}{\gamma}$ may contain fewer elements of $\V_j$ than $\tau$). 
But this may only happen if $w = \Gap(\tau)$ and occurs in $\V_t$ with 
$t < j = i$. Hence $ \Gap(\tau)=w \succ \Gap(\bk{\tau}{\gamma}) \succ u$, 
a contradiction
to our assumption that $u \notin \tail(\tau)$.

\smallskip\noindent{\bf Case 2:} $u \in \up(\bk{\tau}{\gamma})$. 
Let $j$ be the index such that $u \in \V_j$.  
If $j \in \FL(\tau)$, then $u \in \up(\tau)$ 
(as $\fgap(\bk{\tau}{\gamma}, j) \preceq \fgap(\tau, j)$), and the assertion follows. 
So suppose $j \notin \FL(\tau)$. Then $w \in \V_j$, $j < i$ 
(as this is the only way in which $\bk{\tau}{\gamma}$ may contain more elements 
from $\V_j$ than does $\tau$). But then $\fgap(\bk{\tau}{\gamma}, j)$ is not in 
$\tau$: this is because the only vertex that is not in 
$\bk{\tau}{\gamma}$ but is in $\tau$ is $v$, and $v$ is not in $\V_j$. 
Since $\fgap(\bk{\tau}{\gamma}, j) \succ u$ and $j \notin \FL(\tau)$, 
we obtain that $u \in \tail(\tau)$.
\endproof

\begin{lemma} \label{R-back2} 
Let $v$ and $w$ be as in the definition of $\bk{\tau}{\gamma}$,
and suppose $v \in \V_i$, $w \in \V_j$. Then for $1 \leq t \leq m$, $t \neq i$, 
$\fll{\tau}{t}{a_t} = \fll{\bk{\tau}{\gamma}}{t}{a_t}$, 
and so in particular $\mu_{\X_t} = \mu'_{\X_t}$. 
Furthermore $\fll{\tau}{i}{a_i} \preceq \fll{\bk{\tau}{\gamma}}{i}{a_i}$, 
and $\deg(\mu'_{\X_i}) \leq \deg(\mu_{\X_i})$.
\end{lemma}
\proof  If $t \neq i, j$, then $\tau \cap \V_t = \bk{\tau}{\gamma} \cap \V_t$, 
and hence $\fll{\tau}{t}{a_t} = \fll{\bk{\tau}{\gamma}}{t}{a_t}$. 

If $t = j \neq i$, we again have that 
$\fll{\tau}{t}{a_t} = \fll{\bk{\tau}{\gamma}}{t}{a_t}$: indeed, 
as $\bk{\tau}{\gamma} \cap \V_j$ is obtained from $\tau \cap \V_j$ 
by adding the first missing element, 
the first $a_j$-subsets of $\V_j$ containing them are the same.

Now consider $\fll{\bk{\tau}{\gamma}}{i}{a_i}$ and $\fll{\tau}{i}{a_i}$. 
If $v$ is in the initial segment of $\fll{\tau}{i}{a_i}$ 
then $\fll{\bk{\tau}{\gamma}}{i}{a_i}= \fll{\tau}{i}{a_i}$. 
Otherwise, $\fll{\bk{\tau}{\gamma}}{i}{a_i}$ is obtained from 
$\fll{\tau}{i}{a_i}$ by removing $v$ and replacing it with the 
first vertex in $\V_i$ missing from $\fll{\tau}{i}{a_i}$. 
Then $\fll{\tau}{i}{a_i} \preceq \fll{\bk{\tau}{\gamma}}{i}{a_i}$ 
and it follows from the $m = 1$, $\X_0 = \emptyset$ case of 
Lemma \ref{R-back} that $\deg(\mu'_{\X_i}) \leq \deg(\mu_{\X_i})$.
\endproof

\begin{corollary} \label{last-cor}
With the notation as in the previous lemma, 
$\deg(\mu'_{(\X_0 \cup \X_i)}) \leq  \deg(\mu_{(\X_0 \cup \X_i)})$. 
Furthermore, there is a divisor $\nu$ of $\mu_{(\X_0 \cup \X_i)}$ 
with degree equal to $\deg(\mu'_{(\X_0 \cup \X_i)})$, such that 
$\mu'_{(\X_0 \cup \X_i)} \succeq \nu$.
\end{corollary}
\proof The first claim follows from Lemma \ref{R-back} 
and the first part of Lemma \ref{R-back2}. 

To see the second, we note that by the second part of 
Lemma \ref{R-back2} and Lemma \ref{littlephi}, there is a divisor 
$\nu_i$ of $\mu_{\X_i}$ with degree equal to that of $\mu'_{\X_i}$ 
and $\mu'_{\X_i} \succeq \nu_i$. Now let $\nu$ be the revlex last 
divisor of $\mu_{(\X_0 \cup \X_i)} $ with degree equal to that of 
$\mu'_{(\X_0 \cup \X_i)}$. With $\nu_i$ chosen as in the proof of 
Lemma \ref{littlephi}, it is clear that $\nu_i$ divides $\nu$, 
and as all the variables in $\X_i$ come after those in $\X_0$, 
we must have $\mu'_{(\X_0 \cup \X_i)} \succeq \nu$.
\endproof

Lemmas \ref{inverse}--\ref{R-back2} and Corollary \ref{last-cor} put
together imply that $\Phi$ satisfies 
Theorem \ref{bijection} with $\tau' = \bk{\tau}{\gamma}$ 
and $\mu'' = \nu \prod_{1 \leq t \leq m, t \neq i} \mu_{\X_t}$.
This completes the proof of Theorem \ref{main}.


{\small \begin{thebibliography}{999}
\bibitem{AHH}  A.~Aramova, J.~Herzog, and T.~Hibi, 
Shifting operations and graded Betti numbers,
J.~Algebraic Combin.~12 (2000), 207--222.

\bibitem{BjFrSt}
A.~Bj{\"o}rner, P.~Frankl, and R.~Stanley,
The number of faces of balanced Cohen-Macaulay complexes 
and a generalized Macaulay theorem, Combinatorica
7 (1987), 23--34.

\bibitem{Br}
J.~Browder, Shellable complexes from multicomplexes, 
 arXiv:0812.4562. 
 
\bibitem{ClLi}
G.~F. Clements and B.~Lindstr{\"o}m,
A generalization of a combinatorial theorem of {M}acaulay,
 J. Combinatorial Theory 7 (1969), 230--238.


\bibitem{Eis}
D.~Eisenbud,
Commutative Algebra with a View Toward 
Algebraic Geometry, Springer, New York, 1995.

\bibitem{FrFuKa}
P.~Frankl, Z.~F\"uredi, and G.~Kalai, 
Shadows of colored complexes,
Math.~Scand.~63 (1988), 169--178. 


\bibitem{Kat}
G.~O.~H.~Katona, A theorem of finite sets, in Theory
of graphs (Proc.~Colloq.~Tihany, 1966), 
Academic Press, New York, 1968.

\bibitem{Krus}
J.~Kruskal, The number of simplices in a complex,
Mathematical Optimization Techniques, University
of California Press, Berkeley and Los Angeles, 1963,
pp.~251--278.

\bibitem{Macaulay}
F.~S.~Macaulay, Some properties of enumeration in the
theory of modular systems, Proc.~London Math.~Soc.~26 (1927), 531--555.


\bibitem{Mc}
P.~McMullen, The maximum number of faces of a convex polytope,
Mathematika 17 (1970), 179--184.

\bibitem{MePe}
J.~Mermin and I.~Peeva, Lexifying ideals,
Math. Res. Lett., 13(2-3) (2006), 409--422.

\bibitem{Nov05}
I.~Novik, On face numbers of manifolds with symmetry,
Adv.~Math.~192 (2005), 183--208. 

\bibitem{Reis}
G.~A.~Reisner, Cohen-Macaulay quotiens of polynomial rings,
Adv.~Math.~21 (1976), 30--49.


\bibitem{St77}
R.~Stanley, Cohen-Macaulay complexes, in: M.Aigner (Ed.),
Higher Combinatorics, 
Reidel, Dordrecht and Boston, 1977, pp.~51--62.

\bibitem{St79}
R.~Stanley, Balanced Cohen-Macaulay complexes, 
Trans.~Amer.~Math.~Soc. 249 (1979), 139--157.

\bibitem{St96}
R.~P.~Stanley, Combinatorics and Commutative Algebra,
Progress in Mathematics, 41, Birkh{\"{a}}user Boston, Inc., 
Boston, MA, 1996.


\end{thebibliography}}
\end{document}